\documentclass[11pt]{article}
\usepackage{epsfig}
\usepackage{amssymb}
\usepackage{amscd}
\usepackage{epsf}
\usepackage{amsfonts,amssymb,amsthm,amsmath}

\long\def\comment#1\endcomment{\relax}

\newcounter{subsubsubsection}
\newcounter{subsubsubsubsection}

\makeatletter
\@addtoreset{subsubsubsection}{subsubsection}
\@addtoreset{subsubsubsubsection}{subsubsubsection}

\makeatother

\newcommand{\sevafigc}[4]{\begin{figure}[h]\centerline{
 \epsfig{file=#1,width=#2,angle=#3}}
\bigskip\caption{#4}\end{figure}}

\DeclareMathSizes{11.1}{10}{8}{6}

\DeclareMathOperator{\Hom}{Hom}

\newtheorem*{theorem*}{Theorem}

\newtheorem*{lemma}{Lemma}

\theoremstyle{remark}

\newtheorem*{remark}{Remark}

\newtheorem*{example}{Example}

\theoremstyle{definition}

\newtheorem*{definition}{Definition}

\DeclareMathOperator{\Id}{Id}

\newcommand{\mb}{{\bullet}}

\newcommand\End{\mathrm{End}}
\newcommand\Bi{\mathrm{Bialg}}

\newcommand{\Conf}{\mathrm{Conf}}

\newcommand{\g}{\mathfrak{g}}

\newcommand{\K}{\mathrm{K}}

\newcommand{\St}{\mathrm{St}}

\title{{\tt {\huge A concept of $\frac23$PROP and deformation theory of (co)associative bialgebras}}}
\author{{\tt {\LARGE Boris Shoikhet}}}
\date{November 20, 2003}

\begin{document}\maketitle

\hbox to\textwidth{\hfil\parbox{55mm}{\em To my teacher Borya Feigin in the occasion of his 50th birthday}}
\begin{abstract}
{\tt We introduce a concept of $\frac23$PROP generalizing the Kontsevich concept of $\frac12$PROP.
We prove that some Stasheff-type compactification of the Kontsevich spaces $K(m,n)$
defines a topological $\frac23$PROP structure. The corresponding chain complex is
a minimal model for its cohomology (both are considered as $\frac23$PROPs). We construct a $\frac23$PROP
$\End(V)$ for a vector space $V$. Finally, we construct a dg Lie algebra controlling the deformations
of a (co)associative bialgebra. Philosophically, this construction is a version of the Markl's operadic\
construction from [M1] applied to minimal models of $\frac23$PROPs. }
\end{abstract}
\section*{\tt{Introduction}}
The goal of this paper is to construct a deformation theory for
(co)associative bialgebras. According to general principles, it means that
we are looking for a dg Lie algebra (or, more generally, for an $L_\infty$
algebra) controlling the deformation theory of a (co)associative bialgebra.
In the case of the deformation theory of associative algebras, such a dg Lie
algebra controlling the deformations of an associative algebra $A$, is the
cohomological Hochschild complex of $A$ with the Gerstenhaber bracket. (More
precisely, this Hochschild complex controls the deformations of the
category of $A$-modules).

First of all, recall that a (co)associative bialgebra is a vector space $A$
equipped with the maps $\star\colon A^{\otimes 2}\to A$ (the product) and
$\Delta\colon A\to A^{\otimes 2}$ (the coproduct). The product is supposed
to be associative and the coproduct is supposed to be coassociative.
Moreover, we suppose the following compatibility of them:
\begin{equation}\label{eq230.1}
\Delta(a\star b)=\Delta(a)\star\Delta(b)
\end{equation}
for any $a,b\in A$. (Here in the r.h.s. the product is the component product
in $A^{\otimes 2}$ defined as $(a\otimes b)\star(a_1\otimes b_1)=
(a\star a_1)\otimes (b\star b_1)$). Notice that we do not suppose the
existence of unit and counit in $A$.

Here we meet our first difficulty: the r.h.s. of (\ref{eq230.1}) is of the
4th degree and not quadratic. Recall that (little bit roughly) we associate
the deformation theory with a dg Lie algebra $\g^\mb$ as follows: we
consider the solutions of the Maurer-Cartan equation
\begin{equation}\label{eq230.2}
d\alpha+\frac12[\alpha,\alpha]=0
\end{equation}
for $\alpha\in \g^1$ modulo the action of the gauge group associated with
$\g^0$ on these solutions. (Because of possible divergences in the action of
the gauge group, we say instead of this direct construction that the
deformation functor is a functor from the category of the Artinian algebras
to the category of sets).

It is known that the Gerstenhaber-Schack complex [GS] associated with a
bialgebra $A$ is a deformation complex of the bialgebra structure on $A$.
It means that the first cohomology of this complex are isomorphic to the
infinitesimal deformations on $A$. To pass from the infinitesimal
deformations to the global ones, one needs to have an appropriate dg Lie algebra
structure on the Gerstenhaber-Schack complex (or, more generally, an
$L_\infty$-structure). Recall here that as a vector space, the
Gerstenhaber-Schack complex of $A$ is
\begin{equation}\label{eq230.3}
K^\mb_{GS}=\bigoplus_{m,n\ge 1}\Hom(A^{\otimes m},A^{\otimes n})[-m-n+2]
\end{equation}

In particular, in degree 1 we have: $K^1_{GS}=\Hom(A^{\otimes
2},A)\oplus\Hom(A,A^{\otimes 2})$. We could expect that for some dg Lie
algebra structure on $K^\mb_{GS}$ the Maurer-Cartan equation (\ref{eq230.2})
for the element $*_1\oplus \Delta_1\in\Hom(A^{\otimes
2},A)\oplus\Hom(A,A^{\otimes 2})$ means exactly that
$(*+*_1,\Delta+\Delta_1)$ defines a new (co)associative bialgebra structure
on $A$.

But it is impossible: because the r.h.s. of the equation (\ref{eq230.1}) is
of the 4th degree in $*_1$ and $\Delta_1$, while the Maurer-Cartan equation
(\ref{eq230.2}) is quadratic. It means that the best we could expect is to
have an $L_\infty$ algebra structure on $K^\mb_{GS}$ (which looks quite
complicated). This crucial observation was explained to the author by Boris
Tsygan about 3 years ago.

Now remember that the $L_\infty$ algebras and the dg Lie algebras is more or
less the same: if we have an $L_\infty$ algebra structure on a graded vector
space $V$, we necesserily have an $L_\infty$ isomorphic structure of pure dg
Lie algebra on a (bigger) space $V_1$. It means that the question of which
structure we have, dg Lie algebra or $L_\infty$ algebra, is the question of
the right choice of "generators". This means, in particular, that we could
expect the existence of a complex quasi-isomorphic to the
Gerstenhaber-Schack complex ("with another generators") and a dg Lie algebra
on it, which solves the deformation problem of a (co)associative bialgebra
$A$.

This idea is one source of the theory developed here in this paper. Another
source is the Kontsevich spaces $K(m,n)$. The reader can find the definition
of them in Section 2 of the paper. The original Kontsevich motivation when
he invented these spaces was the following:
the space $K(2,2)$ is the configuration space on two independent lines, we
have 2 points in each line modulo independent common shift on each line, and
modulo the following action of $\mathbb{R}^*_+$ on this space: for
$\lambda\in\mathbb{R}^*_+$, we dilatate the first line with the scale
$\lambda$ and the second with the coefficient $\lambda^{-1}$. Then $K(2,2)$
is a 1-dimensional space: we have an interval on the first line, an interval
on the second (we identify the intervals of the same length), and we
identify such two configurations with the same product of the lengthes of
the two intervals. Therefore, the configuration has the only one module--the
product of the lengthes of the intervals. Before compactification, it is
isomorphic to $\mathbb{R}_+$.

Now we compactify the space $K(2,2)$ to the closed interval. The two limit
configurations are shown in the Figure 1:
\sevafigc{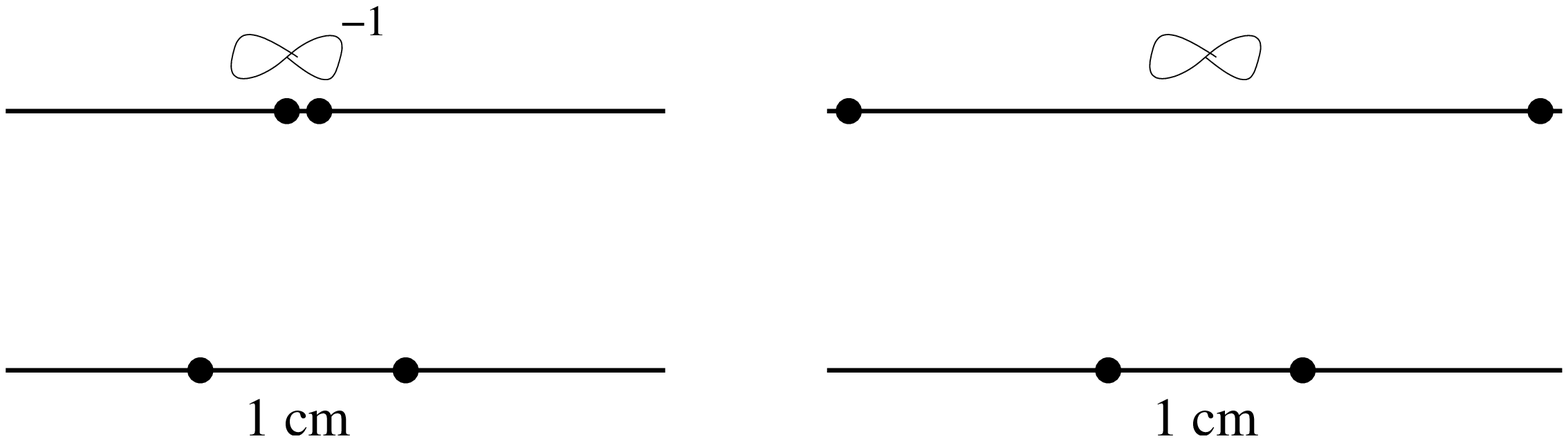}{100mm}{0}{The two limit points in $\overline{K(2,2)}$}
The Kontsevich's insight was that the left picture should give the left-hand
side of the compatibility equation (\ref{eq230.1}), while the right picture
should give the right-hand side of (\ref{eq230.1}). We say "should give"
having in mind the Markl's construction in [M1], or even further, a
construction of the type of Kontsevich formality.

After these short remarks, we pass to our constructions.

If the reader is interested mostly in our construction of the deformation dg
Lie algebra, he can begin to read the paper from Section 3 and to come back
to the previous two Sections if it is necesserily.
\section{\tt{The concept of $\frac23$PROP}}
\subsection{{\tt The definition}}
Here we define our main technical tool--$\frac23$PROPs. The name
$\frac23$PROPs indicates that this concept is a further generalization (or
simplification) of the concept of $\frac12$PROP due to Maxim Kontsevich
(see [K3], [MV]).

We define a $\frac23$PROP of vector spaces, a $\frac23$PROP of dg vector
spaces, of topological spaces,... can be defined analogously.
\begin{definition}
A pre-$\frac23$PROP of vector spaces consists of the following data:
\begin{itemize}
\item[(i)] a collection of vector spaces $F(m,n)$ defined for
$m\ge1,n\ge1,m+n\ge3$, with an action of symmetric groups
$\Sigma_m^\vee\times\Sigma_n$ on $F(m,n)$,
\item[(ii)] a collection of vector spaces $F^{1,1,\dots,1\ (n\ times)}_m$
defined for $n\ge 2,m\ge 2$, with an action of the symmetric group
$\Sigma_m^\vee$ on $F^{1,1,\dots,1\ (n\ times)}_m$,
\item[(iii)] a collection of vector spaces $F_{1,1,\dots,1\ (m\ times)}^n$
defined for $m\ge 2,n\ge 2$, with an action of the symmetric group
$\Sigma_n$ on $F_{1,1,\dots,1\ (m\ times)}^n$,
\item[(iv)] compositions $\circ_i\colon F(m,n)\otimes F(1,n_1)\to
F(m,n+n_1-1)$, $1\le i\le n$,
\item[(v)] compositions $_j\circ\colon F(m_1,1)\otimes F(m,n)\to
F(m+m_1-1,n)$, $1\le j\le m$,
\item[(vi)] compositions $\circledcirc_i\colon F_{1,1,\dots,1\ (m\
times)}^n\otimes F_{1,1,\dots,1\ (m\
times)}^{n_1}\to F_{1,1,\dots,1\ (m\
times)}^{n+n_1-1}$, $1\le i\le n$,
\item[(vii)] compositions $_j\circledcirc\colon F^{1,1,\dots,1\ (n\
times)}_{m_1}\otimes F^{1,1,\dots,1\ (n\
times)}_{m}\to F^{1,1,\dots,1\ (n\
times)}_{m+m_1-1}$, $1\le j\le m$,
\item[(viii)] compositions $\circledcirc\colon F^{1,1,\dots,1\ (n\
times)}_{m}\otimes F_{1,1,\dots,1\ (m\
times)}^n\to F(m,n)$,
\item[(ix)] all the compositions are equivariant with respect to the actions
of the symmetric groups.
\end{itemize}
This data should obey the following properties:
\begin{itemize}
\item[(1)] the composition $\circ_j\star\circ_i\colon F(m,n)\otimes
F(1,n_1)\otimes F(1,n_2)\to F(n,n+n_1+n_2-2)$ is associative for
$i\le j\le i+n_1$: in this case $\circ_{j-i}\star\circ_i\colon
F(m,n)\otimes
(F(1,n_1)\otimes F(1,n_2))\to F(n,n+n_1+n_2-2)$ coincides with
$\circ_j\star\circ_i \colon  (F(m,n)\otimes
F(1,n_1))\otimes F(1,n_2)\to F(n,n+n_1+n_2-2)$,
\item[(2)] in the notations of (1), if $j<i$ or $j>i+n_1$, we have
$\circ_j\star\circ_i\colon (F(m,n)\otimes
F(1,n_1))\otimes F(1,n_2)\to F(n,n+n_1+n_2-2)$ is equal to
$\circ_{i_1}\star\circ_{j_1}\colon (F(m,n)\otimes
F(1,n_2))\otimes F(1,n_1)\to F(n,n+n_1+n_2-2)$ where
$j_1=j, i_1=i+n_2$ if $j<i$, and $j_1=j-n_1, i_1=i$ if $j>i+n_1$ {\it (the commutativity)},
\item[(3)] the property analogous to (1) for $_j\circ$,
\item[(4)] the property analogous to (2) for $_j\circ$,
\item[(5)-(8)] the analogous proprties for $\circledcirc_i$ and for
$_j\circledcirc$.
\end{itemize}
\end{definition}

Notice that this structure without $F^{1,1,\dots,1}_m$ and
$F_{1,1,\dots,1}^n$ is exactly the Kontsevich's $\frac12$PROP structure.
\begin{definition}
A pre-$\frac23$PROP is a $\frac23$PROP if the operations of $\circ$-type are
compatible with the operations of $\circledcirc$-type, as follows:

There are extra maps $\curlywedge_{m_1\to m_2}^n\colon F_{1,1,\dots,1\ (m_1\
times)}^n\to F_{1,1,\dots,1\ (m_2\ times)}^n$, $(m_1\le m_2)$ and mapd
$\curlyvee_m^{n_1\to n_2}\colon F_m^{1,1,\dots,1\ (n_1\ times)}\to F_m^{1,1,\dots,1\ (n_2\
times)}$, $n_1\le n_2$, which are supposed to be equivariant with respect to
the actions of symmetric groups. We also suppose that these maps are isomorphisms. Then we have:
\begin{itemize}
\item[(A)]$\circ_i \star\circledcirc\colon(F^{1,1,\dots,1\ (n\ times)}_m\otimes F_{1,1,\dots,1\
(m\ times)}^n)\otimes F(1,n_1)\to F(m,n+n_1-1)$ is equal to
$\circledcirc\star\circledcirc_i\colon F_m^{1,1,\dots,1\ (n+n_1-1 \ times)}\otimes (F_{1,1,\dots,1\ (m\ times)}^n\otimes
F_{1,1,\dots,1\ (m\ times)}^{n_1})\to F(m,n+n_1-1)$. More precisely, let $\alpha\in F^{1,1,\dots,1\ (n\
times)}_m$, $\beta\in F_{1,1,\dots,1\
(m\ times)}^n)$, and $\gamma\in F(1,n_1)$. Then
\begin{equation}\label{eq23.01}
(\alpha\circledcirc\beta)\circ_i\gamma=  \curlyvee_m^{n\to
n+n_1-1}(\alpha)\circledcirc_i(\beta\circledcirc \curlywedge^{n_1}_{1\to
m}(\gamma))
\end{equation}
\item[(B)] the analogous compatibility with $_j\circ$.
\end{itemize}
\end{definition}
We can imagine what is a free $\frac23$PROP. It consists from all "free"
words of the following two forms:
\begin{equation}\label{eq23.1}
\dots\circ F(m_k,1)\circ\dots\circ F(m_1,1)\circ F(m,n)\circ
F(1,n_1)\circ\dots\circ F(1,n_\ell)\circ\dots
\end{equation}
and
\begin{multline}\label{eq23.2}
 F^{1,1,\dots,1\ (n_1+n_2+\dots-k+1 \ times)}_{m_k}\circledcirc
\dots\circledcirc F^{1,1,\dots,1\ (n_1+n_2+\dots-k+1 \
times)}_{m_1}\circledcirc\\
\circledcirc
F_{1,1,\dots,1 \ (m_1+m_2+\dots-l+1 \ times)}^{n_1}\circledcirc
 F_{1,1,\dots,1 \ (m_1+m_2+\dots -l+1\
times)}^{n_2}\circledcirc\dots\circledcirc F_{1,1,\dots,1 \ (m_1+m_2+\dots \ times)}^{n_l}
\end{multline}

We can draw these free elements as "two-sided trees", see Figure 2:
\sevafigc{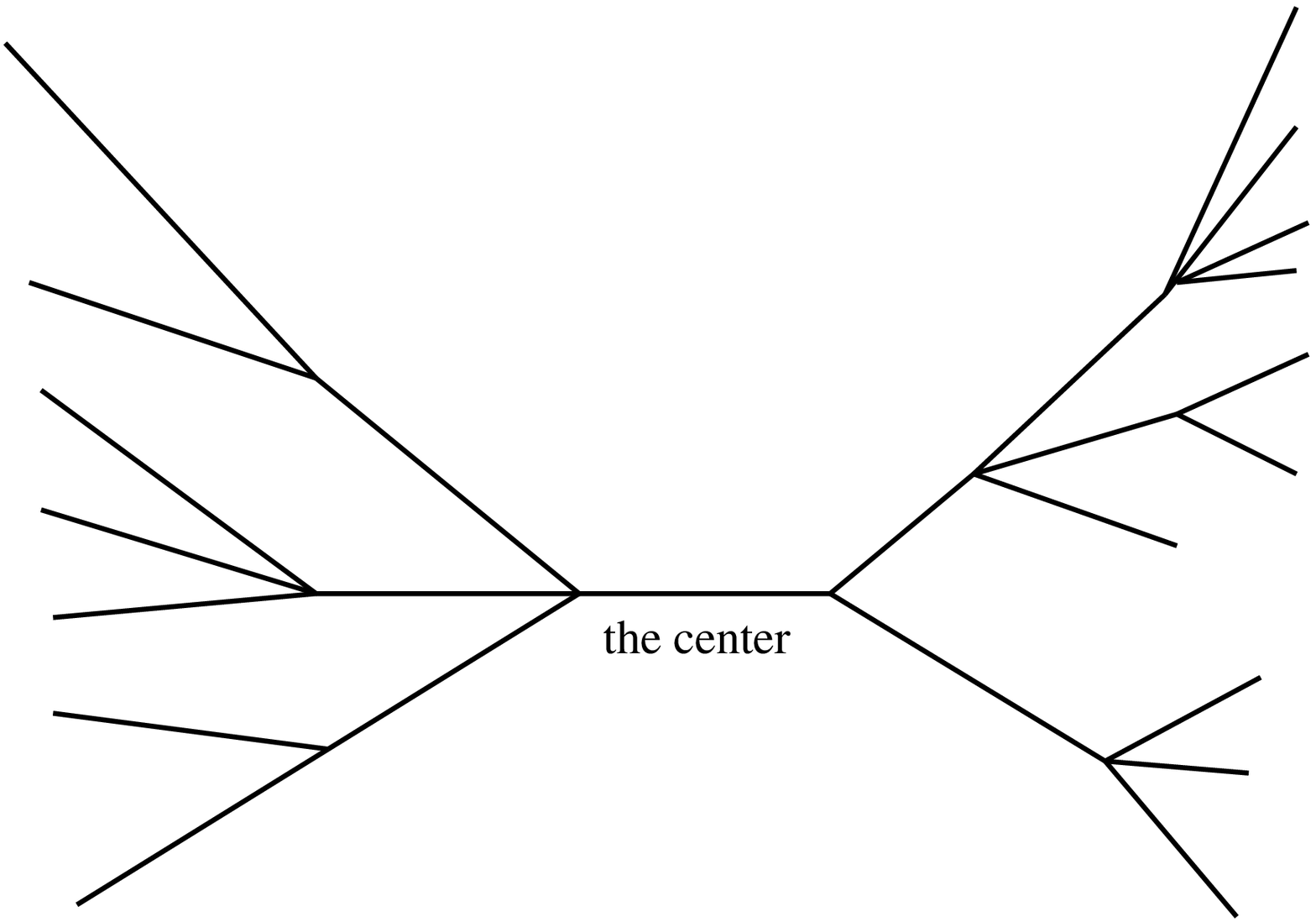}{90mm}{0}{A typical element of a free $\frac23$PROP}
The valences of the vertices are greater or equal than 3. We have trees of
"two different colors" according to the equations
(\ref{eq23.1})-(\ref{eq23.2}) above.

The definition of (pre-)$\frac23$PROP is motivated by the geometry of the
Kontsevich spaces $K(m,n)$ (see the next Section). The reader who is
interested in the origin of this definition can pass directly to Section 2.
We tried to construct a compactification of these spaces and to formalize
the operations among the strata. The advantage of $\frac23$PROPs is that the
chain complex of the compactification $\overline{K(m,n)}$ is a free dg
$\frac23$PROP. Its homology $\frac23$PROP is exactly the $\frac23$PROP $\Bi$
controlling the (co)associative bialgebras. Finally, any (co)associative
bialgebra structure on a vector space $V$ gives a map of the
pre-$\frac23$PROPs $\Bi\to\End(V)$.

\subsection{{\tt The pre-$\frac23$PROP $\End(V)$}}
Here we define the pre-$\frac23$PROP $\End(V)$ for a vector spaces $V$. We use
here the notations $\End(m,n)$, $\End_{1,1,\dots,1}^n$, and
$\End^{1,1,\dots,1}_m$.

We set:
\begin{equation}\label{eq23.3}
\begin{aligned}
\ & \End(m,n)=\Hom(V^{\otimes m},V^{\otimes n}),\\
& \End_{1,1,\dots,1\ (m\ times)}^n=(\Hom(V, V^{\otimes n}))^{\otimes m},\\
& \End^{1,1,\dots,1\ (n\ times)}_m=(\Hom(V^{\otimes m},V))^{\otimes n}
\end{aligned}
\end{equation}
We should now define the compositions $\circ_i,\ _j\circ,\
\circledcirc_i,\ _j\circledcirc$, and $\circledcirc$.

The case of the
composition $\circ_i\colon \End(m,n)\otimes\End(1,n_1)\to\End(m,n+n_1-1)$
can be schematically shown as follows (see Figure 3):
\sevafigc{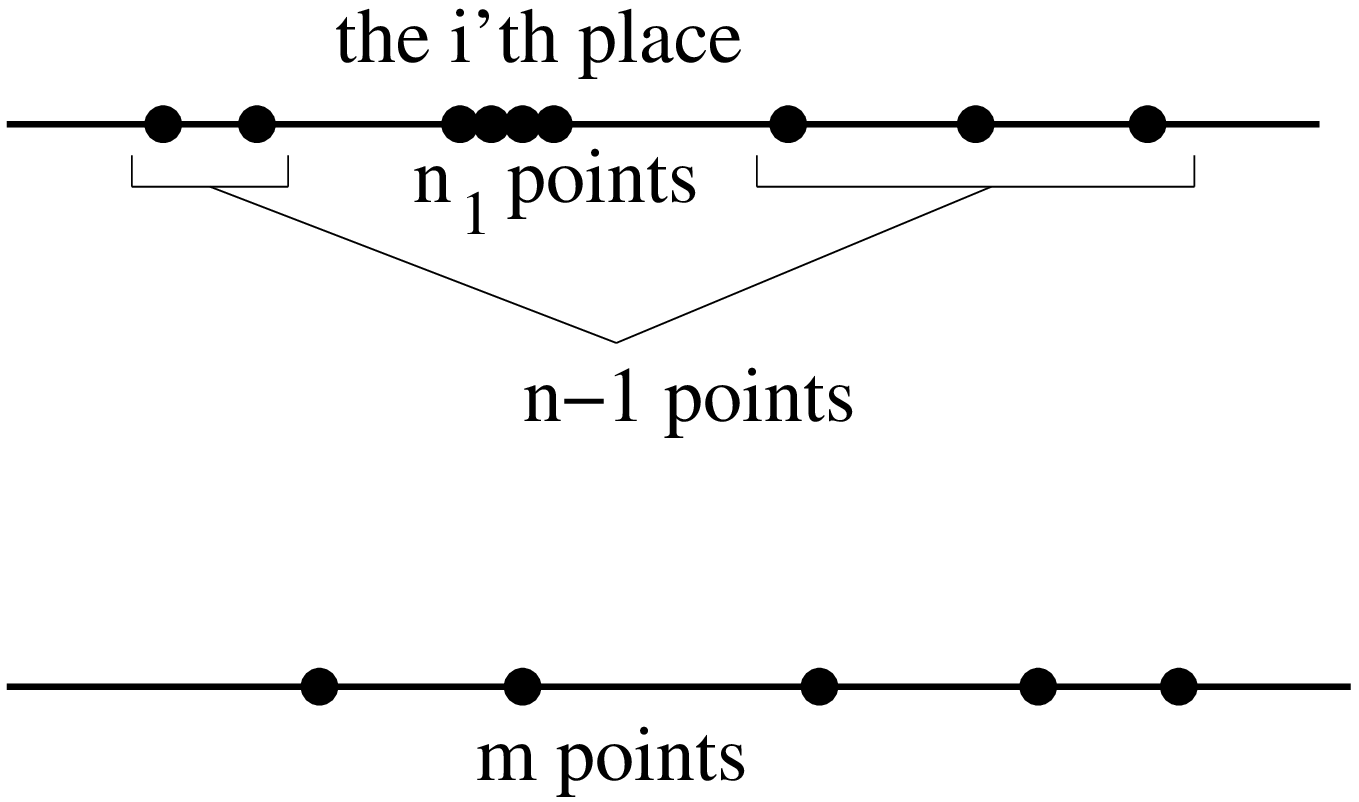}{80mm}{0}{The composition $\circ_i$}
Let $\Psi\in\Hom(V^{\otimes m},V^{\otimes n})$, $\Theta\in\Hom(V,V^{\otimes
n_1})$. Their composition $\Psi\circ_i\Theta\in\Hom(V^{\otimes m},
V^{\otimes n+n_1-1})$ is defined as
\begin{equation}\label{eq232.1}
\Psi\circ_i\Theta (v_1\otimes\dots\otimes v_m)=
(\Id\otimes\dots\otimes\Id\otimes\Theta\otimes\dots\otimes\Id)\circ\Psi(v_1\otimes\dots\otimes
v_m)
\end{equation}
(Here $\Theta$ stands at the $i$th place).

The picture for the composition
$_j\circ\End(m_1,1)\otimes\End(m,n)\to\End(m+m_1-1,n)$ is the following (see
Figure 4):
\sevafigc{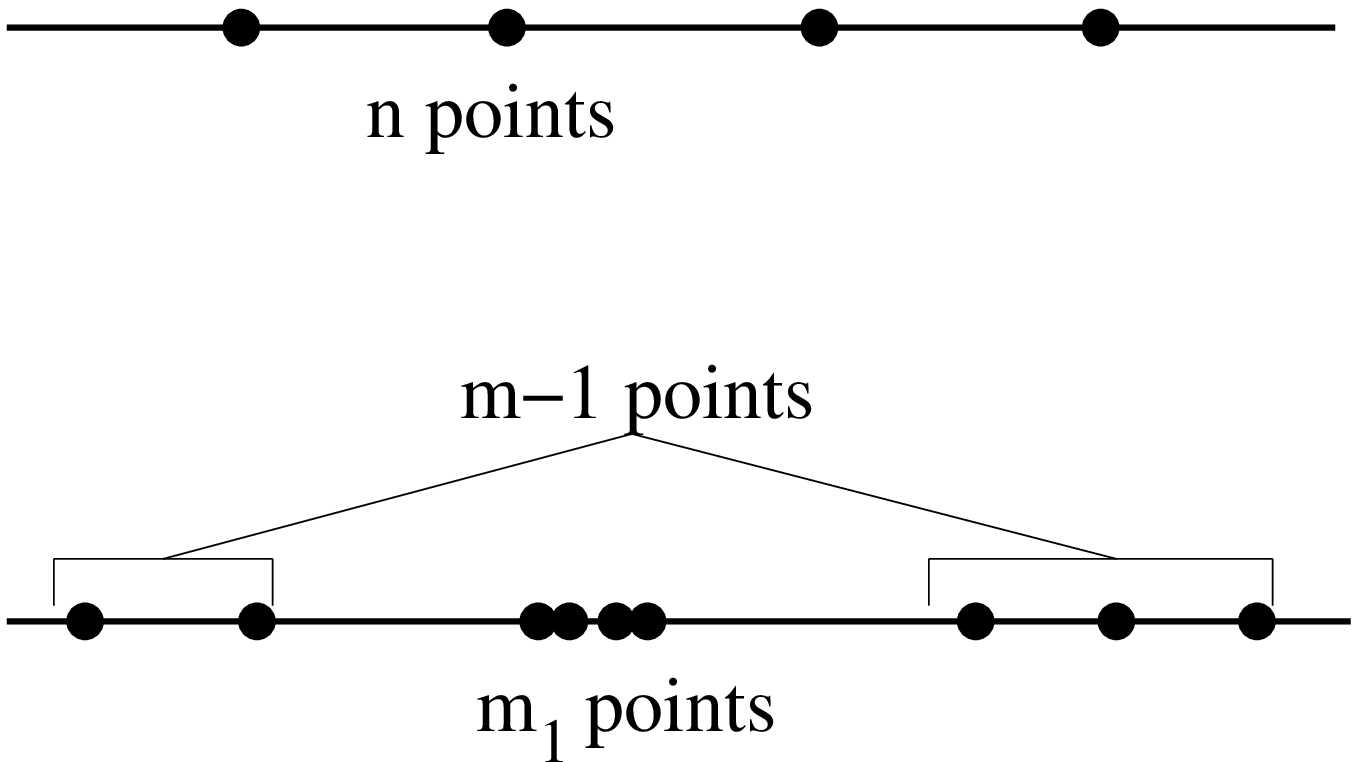}{80mm}{0}{The composition $_j\circ$}
For $\Theta\in \Hom(V^{\otimes m_1}, V)$ and $\Psi\in\Hom(V^{\otimes
m},V^{\otimes n})$, their composition
$\Theta(_j\circ)\Psi\in\Hom(V^{\otimes m+m_1-1},V^{\otimes n})$ is
\begin{multline}\label{eq232.2}
\Theta(_j\circ)\Psi(v_1\otimes\dots\otimes
v_{m+m_1-1})=\\
\Psi(v_1\otimes\dots\otimes
v_{j-1}\otimes\Theta(v_j\otimes\dots\otimes v_{j+m_1-1})\otimes
v_{j+m_1}\otimes\dots\otimes v_{m+m_1-1})
\end{multline}

We have the particular case of the product $\circ_i$ when $m=1$. Denote it
by $\circ_i^1$. By definition, the composition
$\circledcirc_i\colon\End_{1,1,\dots,1\ (m\ times}^n\otimes \End_{1,1,\dots,1\ (m\
times}^{n_1}\to \End_{1,1,\dots,1\ (m\ times}^{n+n_1-1}$ is the $m$th tensor
pover of the composition $\circ_i^1$. Analogously we define $_j\circ^1$ and
the composition $_j\circledcirc\colon \End_{m_1}^{1,1,\dots,1\ (n\
times)}\otimes \End_{m}^{1,1,\dots,1\ (n\
times)}\to \End_{m_1+m-1}^{1,1,\dots,1\ (n\
times)}$ as the $n$th tensor power of the composition $_j\circ^1$.

It remains to define the composition
$\circledcirc$. We define the composition $\circledcirc\colon
(\Hom(V^{\otimes m},V))^{\otimes n} \otimes (\Hom(V,V^{\otimes n}))^{\otimes
m}\to\Hom(V^{\otimes m},V^{\otimes n})$. Suppose
$\Psi_1\otimes\dots\otimes\Psi_n\in(\Hom(V^{\otimes m},V))^{\otimes n}$, and
$\Theta_1\otimes\dots\otimes\Theta_m\in(\Hom(V,V^{\otimes n}))^{\otimes
m}$ We are going to define their composition
$(\Psi_1\otimes\dots\otimes\Psi_n)\circledcirc(\Theta_1\otimes\dots\otimes\Theta_m)\in\Hom(V^{\otimes
m},V^{\otimes n}$.
Denote
\begin{equation}\label{eq232.3}
F(v_1\otimes\dots\otimes
v_m)=\Theta_1(v_1)\otimes\dots\otimes\Theta_m(v_m)\in V^{\otimes mn}
\end{equation}
Next, we define a map $G\in\Hom(V^{\otimes mn}, V^{\otimes n})$ as follows:
\begin{multline}\label{eq232.4}
G(w^1_1\otimes\dots\otimes w_1^n\otimes w_2^1\otimes\dots\otimes
w_2^n\otimes\dots\otimes w_m^1\otimes\dots\otimes w_m^n)=\\
\Psi_1(w_1^1\otimes w_2^1\otimes
w_m^1)\otimes\Psi_2(w^2_1\otimes\dots\otimes w^2_m)\otimes\dots\otimes
\Psi_n(w^n_1\otimes\dots\otimes w^n_m)\in V^{\otimes n}
\end{multline}
Now we set
\begin{equation}\label{eq232.5}
(\Psi_1\otimes\dots\otimes\Psi_n)\circledcirc(\Theta_1\otimes\dots\otimes\Theta_m)=
(G\circ F)(v_1\otimes\dots\otimes v_m)
\end{equation}
It is clear that these compositions define a pre-$\frac23$PROP structure on
$\End(V)$.
\begin{remark}
M.~Markl communicated to the author that our composition $\circledcirc$ is a
particular case of his "fractions" composition [M2].
\end{remark}
\subsection{\tt{The $\frac23$PROP of (co)associative bialgebras $\Bi$}}
For a (pre-)$\frac23$PROP $F$ we define an $F$-algebra structure on a vector
space $V$ as a map of pre-$\frac23$PROPs $F\to\End(V)$. We are going to
construct now a $\frac23$PROP $\Bi$ such that a $\Bi$-algebra structure on
$V$ is exactly a (co)associative bialgebra structure on $V$.

Let $\Sigma_n$ be the symmetric group on $n$ points, and for a group $G$
denote by $G^\vee$ the dual group.

We can consider $\Bi$ as $\frac23$PROP of sets, or, if we like, as the
corresponding $\frac23$PROP of vector spaces (generated by these sets).
We here consider $\Bi$ as a $\frac23$PROP of sets. Later it will
appear also as the homology $\frac23$PROP of the topological $\frac23$PROP
$\overline{K(m,n)}$, then we consider it as the corresponding $\frac23$PROP
of vector spaces. This functor replaces the direct product $\times$ to the
tensor product $\otimes$.

We set:
\begin{equation}\label{eq232.6}
\begin{aligned}
\ &\Bi(m,n)=\Sigma_m^\vee\times\Sigma_n,\\
& \Bi_{1,1,\dots,1\ (m\ times)}^n=\Sigma_n,\\
&\Bi^{1,1,\dots,1\ (n\ times)}_m=\Sigma_m^\vee
\end{aligned}
\end{equation}
The $\frac23$PROP maps $\curlywedge^n_{m_1\to m}$ and $\curlyvee_m^{n_1\to n}$
are the identity maps. We define the compositions $\circ_i$, $_j\circ$, $\circledcirc_i$, $_j\circledcirc$ and $\circledcirc$ as follows:

Consider any of these compositions for the pre-$\frac23$PROP $\End(V)$, a
composition $\bigstar$.
Suppose that $\Psi\in\End_\alpha$ and $\Theta\in\End_\beta$ are its arguments. These compositions were
defined in the previous Subsection. There is the
$\Sigma_{i_1}^\vee\times\Sigma_{j_1}$-action on $\End_\alpha$ and the
$\Sigma_{i_2}^\vee\times\Sigma_{j_2}$-action on $\End_\beta$. Suppose
$\sigma_1^\vee\times\sigma_1\in\Sigma_{i_1}^\vee\times\Sigma_{j_1}$,
and $\sigma_2^\vee\times\sigma_2\in\Sigma_{i_2}^\vee\times\Sigma_{j_2}$. We
are going to define the composition $(\sigma_1^\vee\times
\sigma_1)\bigstar(\sigma_2^\vee\times\sigma_2)$ in $\Bi$. For this, consider
the composition $((\sigma_1^\vee\times
\sigma_1)\Psi)\bigstar((\sigma_2^\vee\times\sigma_2)\Theta$. It is clear that
it is equal to the action of some $\sigma\in\Bi$ on the product of $\Psi$
and $\Theta$ in $\End$:
\begin{equation}\label{eq232.7}
((\sigma_1^\vee\times
\sigma_1)\Psi)\bigstar((\sigma_2^\vee\times\sigma_2)\Theta=\sigma(\Psi\bigstar\Theta)
\end{equation}
The last equation holds for any $\Psi$ and $\Theta$ in the corresponding
components of $\End$, that is, $\sigma$ does not depend on the choice of
$|psi$ and $\Theta$. We define the composition $(\sigma_1^\vee\times
\sigma_1)\bigstar(\sigma_2^\vee\times\sigma_2)$ as $\sigma$.

It clear that this definition is correct, and in this way we define a
$\frac23$PROP $\Bi$.
\begin{lemma}
A map $\phi\colon \Bi\to\End(V)$ of pre-$\frac23$PROPs is the same that a
(co)associative bialgebra structure on $V$.
\begin{proof}
First, let $V$ be a (co)associative bialgebra with the product $\star$ and
the coproduct $\Delta$. We define a map of pre-$\frac23$PROPs
$\phi_{\star,\Delta}\colon\Bi\to\End(V)$. Let
$\sigma^\vee\times\sigma\in\Bi(m,n)$. We put
$\phi_{\star,\Delta}(\sigma^\vee\times\sigma)(v_1\otimes\dots\otimes v_m)=
\Delta^{n-1}\circ \star^{m-1}(v_1\otimes\dots\otimes v_m)$. Here in the
formula $\Delta^n$ and $\star^m$ are the composition powers of the coproduct
and of the product, correspondingly. Because of the (co)associativity, these
powers are well-defined. Next, for $\sigma\in\Bi_{1,1,\dots,1\
(m\times}^n$ we set $\phi_{\star,\Delta}(\sigma)\in(\Hom(V,V^{\otimes
n}))^{\otimes m}$ is the $m$th tensor power of the map
$v\mapsto \sigma(\Delta^{n-1}(v))$. Analogously, using the product, we
define $\phi_{\star,\Delta}$ on $\Bi^{1,1,\dots,1}_m$. Now we explain why
without the compatibility (\ref{eq230.1}) this definition would be
incorrect.

Consider many identity permutations: $\Id_2\in \Sigma_2=\Bi(1,2)$,
$\Id^2\in \Sigma_2^\vee=\Bi(2,1)$, $\Id_2^{1,1}\in \Bi_{1,1}^2$, and
$\Id^2_{1,1}\in\Bi^{1,1}_2$. Than we have the following identity in $\Bi$:
\begin{equation}\label{eq232.8}
\Id^2\circ\Id_2=\Id^2_{1,1}\circledcirc\Id_2^{1,1}=\Id^\vee\times\Id\in\Bi(2,2)
\end{equation}
It is clear that this identity in $\Bi$ follow some identity in the images
of these elements $\Id_2,\ \Id^2,\ \Id_2^{1,1},\ \Id^2_{1,1}$ by the map
$\phi_{\star,\Delta}$ of pre-$\frac23$PROPs. The reader can easily verify
that this identity is exactly the compatibility (\ref{eq230.1}) in a
(co)associative bialgebra. One can prove also that if the compatibility
holds, the definition of $\phi_{\star,\Delta}$ is correct.

Vice versa, suppose we have a map $\phi\colon\Bi\to\End(V)$ of
pre-$\frac23$PROPs. Denote $a\star b:=\phi(\Id^2)(a\otimes b)$ and
$\Delta(a):=\phi(\Id_2)(a)$. The compatibility follows from (\ref{eq232.8}).
The reader can easily find analogous identities in $\Bi$ which imply the
associativity of $\star$ and the coassociativity of $\Delta$.
\end{proof}
\end{lemma}
\section{\tt{The Kontsevich spaces $K(m,n)$, their Stasheff-type
compactification, and the corresponding $\frac23$PROP}}

First of all, recall the definition of the spaces $\K(m,n)$ due to Maxim
Kontsevich (see also [Sh]). We show in the sequel that these spaces and its
compactification introduced below play a crucial role in the deformation
theory of (co)associative bialgebras.

First define the space $\Conf (m,n)$. By definition, $m,n\ge 1$, $m+n\ge 3$,
and
\begin{multline}\label{eq1.1}
\Conf (m,n)=\{p_1,\dots, p_m\in \mathbb{R}^{(1)}, p_i<p_j\ \  for\ \  i<j;\\
q_1,\dots,q_n\in\mathbb{R}^{(2)}, q_i<q_j\ \ for\ \ i<j\}
\end{multline}
Here we denote by $\mathbb{R}^{(1)}$ and by $\mathbb{R}^{(2)}$ two different
copies of a real line $\mathbb{R}$.

Next, define a 3-dimensional group $G^3$ acting on $\Conf(m,n)$. This group
is a semidirect product $G^3=\mathbb{R}^2\ltimes\mathbb{R}_+$ (here
$\mathbb{R}_+=\{x\in\mathbb{R}, x>0\}$) with the following group law:
\begin{equation}\label{eq1.2}
(a,b,\lambda)\circ(a^{\prime},b^{\prime},\lambda ^{\prime})=
(\lambda ^{\prime} a+a^{\prime},(\lambda ^{\prime})^{-1} b+b^{\prime},\lambda\lambda ^{\prime} )
\end{equation}
where $a,b,a^\prime ,b^\prime\in\mathbb{R}, \lambda,\lambda^\prime\in\mathbb{R}_+$.
This group acts on the space $\Conf(m,n)$ as
\begin{multline}\label{eq1.3}
(a,b,\lambda)\cdot (p_1,\dots,p_m;q_1,\dots,q_n)=
(\lambda p_1+a,\dots,\lambda
p_m+a;\lambda^{-1}q_1+b,\dots,\lambda^{-1}q_n+b)
\end{multline}
In other words, we have two independent shifts on $\mathbb{R}^{(1)}$ and
$\mathbb{R}^{(2)}$ (by $a$ and $b$), and $\mathbb{R}_+$ dilatates
$\mathbb{R}^{(1)}$ by $\lambda$ and dilatates $\mathbb{R}^{(2)}$ by $\lambda
^{-1}$.

In our conditions $m,n\ge 1, m+n\ge 3$, the group $G^3$ acts on $\Conf(m,n)$
freely. Denote by $\K(m,n)$ the quotient-space. It is a smooth manifold of
dimension $m+n-3$.

We will need also a very special case of the spaces  $\K_{m_1,\dots,m_{\ell_1}}^{n_1,\dots,n_{\ell_2}}$  introduced below.
Recall here our definition of the space
$\K_{m_1,\dots,m_{\ell_1}}^{n_1,\dots,n_{\ell_2}}$ (generalizing the
Kontsevich space $\K(m,n)$) from [Sh]:

Fist define the space
$\Conf^{m_1,\dots,m_{\ell_1}}_{n_1,\dots,n_{\ell_2}}$. By definition,
\begin{multline}\label{eq11.4}
\Conf^{m_1,\dots,m_{\ell_1}}_{n_1,\dots,n_{\ell_2}}=\\
\{p^1_1,\dots,p^1_{m_1}\in\mathbb{R}^{(1,1)},
p^2_1,\dots,p^2_{m_2}\in\mathbb{R}^{(1,2)},\dots,
p^{\ell_1}_1,\dots,p^{\ell_1}_{m_{\ell_1}}\in\mathbb{R}^{(1,\ell_1)};\\
q^1_1,\dots,q^1_{n_1}\in\mathbb{R}^{(2,1)},
q^2_1,\dots,q^2_{n_2}\in\mathbb{R}^{(2,2)}\dots,
q^{\ell_2}_1,\dots,q^{\ell_2}_{n_{\ell_2}}\in\mathbb{R}^{(2,\ell_2)}|\\
p^j_{i_1}<p^j_{i_2}\ \ for\ \ i_1<i_2;q^j_{i_1}<q^j_{i_2}\ \ for\ \
i_1<i_2\}
\end{multline}
Here $\mathbb{R}^{(i,j)}$ are copies of the real line $\mathbb{R}$.
Now we have an $\ell_1+\ell_2+1$-dimensional group $G^{\ell_1,\ell_2,1}$
acting on $\Conf^{m_1,\dots,m_{\ell_1}}_{n_1,\dots,n_{\ell_2}}$.
It contains $\ell_1+\ell_2$ independent shifts
$$
p_i^j\mapsto p_i^j+a_j, i=1,\dots, m_j, a_j\in\mathbb{R};
q_i^j\mapsto q_i^j+b_j, i=1,\dots, n_j, b_j\in\mathbb{R}
$$
and {\it one } dilatation
$$
p_i^j\mapsto \lambda\cdot p_i^j\ \ for\ \ all\ \ i,j;q_i^j\mapsto
\lambda^{-1}\cdot q_i^j\ \ for\ \ all\ \ i,j.
$$
This group is isomorphic to $\mathbb{R}^{\ell_1+\ell_2}\ltimes\mathbb{R}_+$.
We say that the lines
$\mathbb{R}^{(1,1)},\mathbb{R}^{(1,2)},\dots,\mathbb{R}^{(1,\ell_1)}$
(corresponding to the factor $\lambda$) are the lines of the first type, and
the lines $\mathbb{R}^{(2,1)},\mathbb{R}^{(2,2)},\dots,\mathbb{R}^{(2,\ell_2)}$
(corresponding to the factor $\lambda^{-1}$) are the lines of the second
type.

Denote
\begin{equation}\label{eq11.5}
\K^{m_1,\dots,m_{\ell_1}}_{n_1,\dots,n_{\ell_2}}=\Conf
^{m_1,\dots,m_{\ell_1}}_{n_1,\dots,n_{\ell_2}}/G^{\ell_1,\ell_2,1}
\end{equation}

We construct a compactication $\overline{K(m,n)}$ the boundary strata of
which are products of the spaces $K(m_1,n_1)$, $K_{1,1,\dots,1}^{n_2}$, and
$K^{1,1,\dots,1}_{m_2}$ (it allowed to be several spaces of each type).
\subsubsection*{{\tt Example}}
Let $m=n=2$. Then the space $\K(2,2)$ is 1-dimensional. It is easy to see
that $(p_2-p_1)\cdot (q_2-q_1)$ is preserved by the action of $G^3$, and it
is the only invariant of the $G^3$-action on $\K(2,2)$. Therefore, $\K(2,2)\simeq
\mathbb{R}_+$. There are two "limit" configurations: $(p_2-p_1)\cdot
(q_2-q_1)\rightarrow 0$ and $(p_2-p_1)\cdot (q_2-q_1)\rightarrow \infty$.
Therefore, the compactification $\overline{\K(2,2)}\simeq [0,1]$. See Figure
1.
\bigskip
We will construct a compactification of the space $K(m,n)$ which is a
topological $\frac23$PROP. More presisely,
\begin{equation}\label{eq2323.1}
\begin{aligned}
\ &F(m,n)=\overline{K(m,n)}\\
&F_{1,1,\dots,1}^n=\overline{K_{1,1,\dots,1}^n}\\
&F_m^{1,1,\dots,1}=\overline{K^{1,1,\dots,1}_m}
\end{aligned}
\end{equation}
The compactifications will be defined below.

First of all, let us describe all strata of codimension 1 in
$\overline{K(m,n)}$. There are boundary strata of codimension 1 of two
different types. The first two strata are shown in Figures 3 and 4. In the
picture in Figure 3 $n_1$ points on the upper line move infinitely close to each other,
and the "scale" of this infinitely small number is irrelevant (we have in
mind here the CROC compactification from [Sh] where this scale is relevant),
$2\le n_1\le n$. In Figure 4 $m_1$ points on the lower line move close to
each other, $2\le m_1\le m$. The remaining stratum of codimension 1 in
$\overline{K(m,n)}$ (there is the only one such stratum) is shown in Figure
5:
\sevafigc{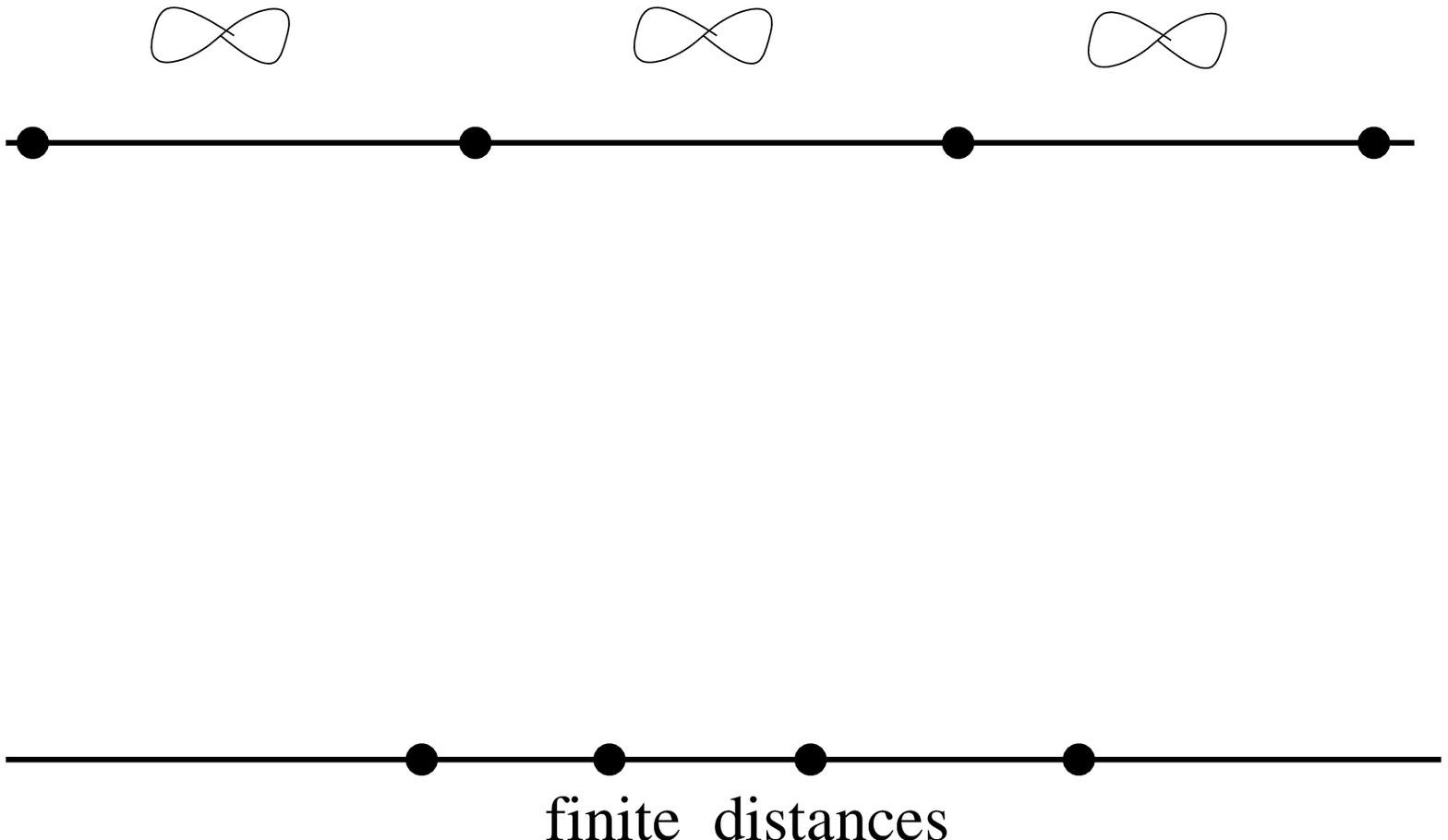}{80mm}{0}{The stratum of codimension 1 in $\overline{K(m,n)}$ of the third type}
Here in the Figure all points on the upper line move infinitely far from
each other {\it will a finite ratio of any two among these infinite
distances}. The all distances on the lower line are finite. Of course, the stratum when
the points on the lower line are infinite and the points on the upper line
are in finite distances from each other, is the same: one stratum can be
obtained from another by the application of the element
$(0,0,\infty)\in G^3$.

The strata in Figures 3 and 4 are isomorphic to $K(m,n)\times K(1,n_1)$ and
$K(m_1,1)\times K(m,n)$, correspondingly. The stratum in Figure 5 is
isomorphic to $K_m^{1,1,\dots,1\ (n\ times)}\times K^n_{1,1,\dots,1\ (m\
times)}$. Let us explain the last formula: the distances between the points
at the upper line are infinite, but their ratios are finite. Therefore, we
should count these ratios. For this, we apply to Figure 5 the transformation
$(0,0,\infty)\in G^3$. Then the infinite distances will become finite, and
then we can count the ratios.

Now we claim that using the 3 operations shown in Figures 3,4,5 we can
obtain {\it any} limit configuration (here by a limit configuration we mean
a configuration where some distances are infinitely large or/and infinitely
small). Moreover, we can apply the configuration in Figure 5 not more than 1
time. Let us explain it:

Apply a transfotmation from $G^3$ such that there are no infinite distances
at the lower line, and the diameter of the configuration of points at the
lower line is finite (not infinitely small). We distinguish the following
two cases: in the first case in the obtained configuration there are no infinitely large distances in the
upper line, and in the second some distances are infinitely large. It is
clear that we can reach any limit configuration of the first type by the
applying several times the degenerations shown in Figure 3 and Figure 4.
Analogously, in the second case, we first apply degenerations in Figure 3
(several times), then apply the transformation from Figure 5 (with the scale
of the infinity depending on the configuration), and then apply several
times the degenerations from Figure 4. It is clear that in this case we can
get {\it any} limit configuration.

Denote the operations shown in Figures 3,4,5 by $\circ_i$, $_j\circ$, and
$\circledcirc$, correspondingly. It remains to define our operations
$\circledcirc_i$ and $_j\circledcirc$ in the Definition in Section 1.1, that
is, to compactify the spaces $K_{1,1,\dots,1}^n$ and $K_m^{1,1,\dots,1}$.

Notice that the space $K^{1,1,\dots,1\ (n\ times)_m}$ is isomorphic to the
Stasheff polyhedron $\St_m$ for {\it any} $n$, as well the space
$K_{1,1,\dots,1\ (m\ times)}^n$ is isomorphic to the Stasheff polyhedron
$\St_n$ for any $m$. In particular, we define our maps $\curlywedge_{m_1\to
m_2}^n$ and $\curlyvee_m^{n_1\to n_2}$ as the identity maps. Furthermore, we
compactify these spaces as usual in the Stasheff compactification, and the
compositions $\circledcirc_i$ and $_j\circledcirc$ are defined in the
natural way. Moreover, it is clear that the formulas (A) and (B) in the
Definition of $\frac23$PROP hold.

Thus, we constructed a topological $\frac23$PROP $\overline{K(m,n)}$.
\begin{lemma}
\begin{itemize}
\item[(i)] The corresponding chain dg $\frac23$PROP (formed by the chain
complexes of the spaces in our stratification) is a free $\frac23$PROP of
graded vector spaces (when we forget about the differential),
\item[(ii)] the corresponding homology $\frac23$PROP is the $\frac23$PROP of
bialgebras $\Bi$
\end{itemize}
\begin{proof}
(i) follows from the fact that any limit configuration can be witten as the
composition of the degenerations in Figures 3,4,5 in a unique way, (ii)
follows from the fact that all spaces $\overline{K(m,n)}$ are contractible
(and, therefore, have only 0-th nontrivial homology).
\end{proof}
\end{lemma}

\begin{example}
In our compactification, the left picture in Figure 1 is $K(2,1)\times
K(1,2)$, and the right picture is $K_{1,1}^2\times K_2^{1,1}$.
\end{example}

Now we are ready to introduce our main object--a dg Lie algebra.

\section{\tt The dg Lie algebra}
Denote by $\mathbb{C}\{\Hom(V^{\otimes m},V)\}$ the vector space generated
by the infinite series of the form
\begin{equation}\label{eq233.1}
\Psi+\Psi\otimes\Psi+\Psi\otimes\Psi\otimes\Psi+\dots\in\prod_{n\ge
1}(\Hom(V^{\otimes m},V))^{\otimes n}
\end{equation}
where $\Psi\in\Hom(V^{\otimes m},V)$ , $m\ge 2$. We denote the sum above by $\overline{\Psi}$.
As a vector space, $\mathbb{C}\{\Hom(V^{\otimes m},V)\}$
is "a very huge" vector space generated by the {\it set} $\Hom(V^{\otimes
m},V)$. Analogously, introduce the notation $\mathbb{C}\{\Hom(V,V^{\otimes
n})\}$, $n\ge 2$, and the notation $\overline{\Theta}$ for $\Theta
\in\Hom(V,V^{\otimes n})$.

Now introduce a dg Lie algebra $\aleph$. First we introduce a graded vector space $\aleph_0$,
and then $\aleph$ will be a quotient space. We set:
\begin{multline}\label{eq233.2}
\aleph_0=\bigoplus_{m\ge 2}\mathbb{C}\{\Hom(V^{\otimes m},V)\}[-m+1]\oplus\bigoplus_{n\ge 2}\mathbb{C}\{\Hom(V,V^{\otimes
n})\}[-n+1]\oplus\\
\oplus\bigoplus_{m,n\ge 2}\Hom(V^{\otimes m},V^{\otimes n})[-m-n+2]
\end{multline}
Now we introduce a prebracket on $\aleph_0$. It means that it is skew-symmetric but does
not obey the Jacobi identity. The idea goes back to the
constructions in Section 2. We associate generators with the strata of
codimension 0, and their compositions (the bracket) by the strata of
codimension 1. Then the Jacobi identity follows (after factorization) from the equation
$\partial^2=0$ where $\partial$ is the chain differential. Philosophically,
it is a kind of the Markl's construction from [M1], but the author can not
verbalize it at the moment. The appearance of the infinite power series
$\overline{\Psi}$ and $\overline{\Theta}$ are motivated by the identity
(\ref{eq23.01}). Thus, the reader can say that no infinite sums appear among
the strata of codimension 0. Nevertheless, the Markl's construction is an
operadic construction, and when we deal with $\frac23$PROPs we need some
modifications.

Only the following brackets are nonzero:
\begin{equation}\label{eq233.3}
\begin{aligned}
\ {\rm (i)}\ \ &
[\overline{\Psi_1},\overline{\Psi_2}]:=\overline{[\Psi_1,\Psi_2]_G}\\
&{\rm where}\ \Psi_1\in\Hom(V^{\otimes m_1},V),\ \Psi_2\in\Hom(V^{\otimes m_2},V), \\
&{\rm and}\ [\Psi_1,\Psi_2]_G\ {\rm is}\ {\rm the}\ {\rm Gerstenhaber}\ {\rm bracket},\\
{\rm (ii)}\ \ &
[\overline{\Theta_1},\overline{\Theta_2}]:=\overline{[\Theta_1,\Theta_2]^G}\\
&{\rm where}\ \Theta_1\in\Hom(V,V^{\otimes n_1}),\ \Theta_2\in\Hom(V,V^{\otimes n_2}),\\
&{\rm and}\ [\Theta_1,\Theta_2]^G\ {\rm is}\ {\rm the}\ {\rm Gerstenhaber}\ {\rm cobracket,}\\
{\rm (iii)}\ \ &
[\overline{\Theta},\overline{\Psi}]:=\Theta\circ\Psi\pm(\Theta)^{\otimes m}\circledcirc
(\Psi)^{\otimes n}\\
&{\rm where}\ \Psi\in \Hom(V^{\otimes m},V),\ \Theta\in\Hom(V,V^{\otimes n}),\
\circ\ {\rm and} \circledcirc\\
&{\rm are}\ {\rm the}\ {\rm compositions}\ {\rm in}\ {\rm the}\ {\rm pre-}\frac23{\rm PROP}\ \End(V),\\
{\rm (iv)}\ \ & {\rm for}\ \alpha\in\Hom(V^{\otimes m},V^{\otimes n}),\ m,n\ge 2,\ \Psi\in\Hom(V^{\otimes m_1},V)\\
&[\overline{\Psi},\alpha]:=\sum_{j=1}^m\pm\Psi(_j\circ)\alpha,\\
{\rm (v)}\ \ & {\rm for}\ \alpha\in\Hom(V^{\otimes m},V^{\otimes n}),\ m,n\ge 2,\ \Theta\in\Hom(V,V^{\otimes n_1})\\
&[\alpha,\overline{\Theta}]:=\sum_{i=1}^n\pm\alpha\circ_i\Theta
\end{aligned}
\end{equation}
We also suppose that the bracket $[\ ,\ ]$ is graded-skew-commutative.

This bracket does {\it not} obey the Jacobi identity. The Jacobi identity
fails for the brackets
$[[\overline{\Psi_1},\overline{\Psi_2}],\overline{\Psi_3}]$,
$[[\overline{\Theta_1},\overline{\Theta_2}],\overline{\Theta_3}]$,
$[[\overline{\Psi_1},\overline{\Psi_2}],\overline{\Theta}]$, and
$[[\overline{\Psi},\overline{\Theta_1}],\overline{\Theta_2}]$.
Our solution is to factorize by the vector space spanned by the Jacobi
identities, as follows:

First, factorize the space generated by $\overline{\Psi}$ by the space $I_G$
{\it generated} by the Jacobi identities. It means that we first consider a
graded vector space $I_G^0$ spanned by the vectors
$\overline{[[\Psi_1,\Psi_2]_G,\Psi_3]_G}\mp\overline{[\Psi_1,[\Psi_2,\Psi_3]_G]_G}\pm
\overline{[\Psi_2,[\Psi_1,\Psi_3]_G]_G}$ where $\Psi_i\in \Hom(V^{\otimes
m_i},V)$ for some $m_1,m_2,m_3\ge 2$. Then consider the vector space $I_G$ spanned
by the elements of the form
$[\overline{\Psi_{k_\ell}},[\overline{\Psi_{k_{\ell-1}}},[\dots[\overline{\Psi_{k_1}},\alpha]]\dots]$
where $\alpha\in I_G^0$, and $k_i\ge 2$ for all $i$, $\ell\ge 0$. Define the
space $\overline{\Hom(V^{\otimes m},V)}$, $m\ge 2$, as the graded component
of the quotient of the space spanned by all vectors $\overline{\Psi}$, by
the graded space $I_G$. It is clear that $\oplus_{m\ge
2}\overline{\Hom(V^{\otimes m},V)}[-m+1]$ with the bracket (\ref{eq233.3}) is a
graded Lie algebra. One can show that this graded Lie algebra is bigger than
$\oplus_{m\ge 2}\Hom(V^{\otimes m},V)[-m+1]$.

Analogously we define the space $I^G$ and the graded Lie algebra $\oplus_{n\ge
2}\overline{\Hom(V,V^{\otimes n})}[-n+1]$.

Now we define the quotient-spaces $\overline{\Hom(V^{\otimes m},V^{\otimes
n})}$ for $m,n\ge 2$.

Notice, that the bracket of $I_G$ and $I^G$ with $\Hom(V^{\otimes
m},V^{\otimes n}),m,n\ge 2$ is {\it zero without any factorization} because it
depends only on the "linear part" of $\overline{\Psi},\overline{\Theta}$.
The same is true for the product $\Psi\circ\Theta$.
Then, define $I_0$ as the graded vector space spanned by all elements of the form
$[\overline{\Psi_1},\overline{\Psi_2}]\circledcirc\overline{\Theta}\mp
\Psi_1\circ_*(\overline{\Psi_2}\circledcirc\overline{\Theta})\pm
\Psi_2\circ_*(\overline{\Psi_1}\circledcirc\overline{\Theta})$ where
$\Psi_i\in\Hom(V^{\otimes m_i},V),m_i\ge 2$, $\Theta\in \Hom(V,V^{\otimes
n}),n\ge 2$ and the analogous expressions with 2 $\Theta$'s and 1 $\Psi$.
Here $A\circ_* B$ denotes the sum with the signs as in (iv),(v) in the
definition (\ref{eq233.3}). Next, denote by $I$ the graded vector space
generated by the elements of the form
$$
\Psi_{k_s}(_*\circ)\Psi_{k_{s-1}}(_*\circ)\dots(_*\circ)\Psi_{k_1}(_*\circ)\alpha\circ_*
\Theta_{\ell_1}\circ_*\dots\circ_*\Theta_{\ell_s},
$$
where $\alpha\in\Hom(V^{\otimes m},V^{\otimes n}),m,n\ge 2$.

Finally, denote by $\overline{\Hom(V^{\otimes m},V^{\otimes n})}[-m-n+2]$, $m,n\ge
2$, the graded component of the quotient of $(\oplus_{m,n\ge
2}\Hom(V^{\otimes m},V^{\otimes n})[-m-n+2])/I$.

Denote
\begin{equation}\label{eq235.1}
\begin{aligned}
\aleph=&\bigoplus_{m\ge 2}\overline{\Hom(V^{\otimes m},V)}[-m+1]\oplus\\
&\bigoplus_{n\ge 2}\overline{\Hom(V,V^{\otimes n})}[-n+1]\oplus\\
&\bigoplus_{m,n\ge 2}\overline{\Hom(V^{\otimes m},V^{\otimes n})}[-m-n+2]
\end{aligned}
\end{equation}

\begin{theorem*}
The formulas for the bracket (\ref{eq233.3}) define a graded Lie algebra structure on $\aleph$
(that means that the Jacobi identity is satisfied).
\begin{proof}
It follows from the definitions.
\end{proof}
\end{theorem*}
\begin{remark}
It would be very interesting to specify in which sense our construction is
an analog of the Markl's construction [M1] applied to the case of
$\frac23$PROPs.
\end{remark}

The dg Lie algebra $\aleph$ defined above is a "deformation Lie algebra"
(see Section 4) for the bialgebra $V$ with 0 product and 0 coproduct. When
we want to consider deformation theory for a bialgebra $V$ with non-zero
(co)product, we localize $\aleph$ by the corresponding solution of the
Maurer-Cartan equation.

\begin{lemma}
Let $\Psi\in\Hom(V^{\otimes 2},V)$ and $\Theta\in \Hom(V,V^{\otimes 2})$ are
the product and the coproduct for a (co)associative bialgebra structure on
$V$. Then $\beta=\overline{\Psi}+\overline{\Theta}\in \aleph^1$ satisfies the
Maurer-Cartan equation with 0 differential:
\begin{equation}\label{eq233.5}
[\beta,\beta]=0
\end{equation}
\begin{proof}
It is clear.
\end{proof}
\end{lemma}

Now for such $\Psi,\Theta$ as above, we consider the dg Lie algebra
$\aleph_{\Psi,\Theta}$ which is the same as $\aleph$ but with the
differential $ad(\beta)$.
\section{{\tt From dg Lie algebra to $L_\infty$ algebra}}
Let $\g_1$, $\g_2$ be two $L_\infty$ algebras. Recall that it means that we
have odd vector fields of degree +1 $Q_1$ on $\g_1[1]$ and $Q_2$ on
$\g_2[1]$ such that $Q_1^2=Q_2^2=0$. Suppose we have an $L_\infty$ map
$U\colon\g_1\to\g_2$. It means, by definition, that we have a non-linear map
$U\colon\g_1[1]\to\g_2[1]$ which maps the field $Q_1$ to the field $Q_2$.
Suppose that the map $U$ is a (non-linear) imbedding of topological spaces. Then we can
say that the vector field $Q_2$ is tangent to the image $U(\g_1[1])$ (because
it coincides with the image of $Q_1$).
Vice versa, suppose we have an $L_\infty$ algebra $\g_2$, a {\it graded
vector space} $\g_1$, and an imbedding $U\colon\g_1[1]\to\g_2[1]$ such that
the vector field $Q_2$ on $\g_2$ is {\it tangent} to the image. Then we
claim that there is a unique $L_\infty$ structure on $\g_1$ which makes $U$
an $L_\infty$ map.

Apply it now to the case when $\g_2=\aleph$. Consider the
Gerstenhaber-Schack space $\bigoplus_{m,n\ge 1, m+n\ge 3}\Hom(V^{\otimes m},
V^{\otimes n})[-m-n+2]$ as $\g_1$. Consider the following non-linear map
$U\colon\g_1\to \aleph$: the map $U$ maps $\Hom(V^{\otimes m},V^{\otimes
n})[-m-n+2]$ identically to $\aleph$ when $m,n\ge 2$. When $n=1$, $U$ maps
$\Psi\in\Hom(V^{\otimes m},V)\in\g_1$ to $\overline{\Psi}\in\aleph$, and for $m=1$,
$U$ maps $\Theta\in\Hom(V,V^{\otimes n})\in\g_1$ to
$\overline{\Theta}\in\aleph$. It is clear that we are in the assumptions
above (before the localization by the solution $\beta$ of the Maurer-Cartan
equation). It means that the quadratic vector field on $\aleph[1]$ defining
the dg Lie algebra structure on $\aleph$ is {\it tangent} to the image of
$U$. It allows us to define an $L_\infty$ structure on the
Gerstenhaber-Schack space $\g_1$ which makes $U$ an $L_\infty$ map.

The author hopes to construct this $L_\infty$ structure explicitly in the
next paper. He is not sure that the linear part of this $L_\infty$ structure
will be the Gerstenhaber-Schack differential, but he is sure that in this
way we obtain a more right object.
\subsection*{Acknowledgements}
I am grateful to Maxim Kontsevich for sharing with me with his construction
of the spaces $K(m,n)$. When I constructed the CROC compactification of
them, I discussed it many times with Borya Feigin, and these discussions
finally explained to me the necessity of a Stasheff-type compactification for
the usual (commutative) deformation theory of (co)associative bialgebras.
Discussions on this compactification and on the corresponding Lie algebra
with Giovanni Felder were very useful
for me. I am grateful to Martin Markl for the reading of the first
version of this paper and for his remarks and corrections. I also would like to thank the ETH (Zurich) for the financial
support and for the very stimulating atmosphere.

\bigskip
\bigskip
Dept. of Math., ETH-Zentrum, 8092 Zurich, SWITZERLAND\\
e-mail: {\tt borya@mccme.ru, borya@math.ethz.ch}

\end{document}